\documentclass[journal]{IEEEtran}

\usepackage{amssymb}
\usepackage[ruled,vlined,linesnumbered]{algorithm2e}
\usepackage{amsthm}
\usepackage{amsfonts}
\usepackage{amsmath}
\usepackage{booktabs}
\usepackage{bm}
\usepackage{color}
\usepackage{cite}
\usepackage{comment}
\usepackage{epsfig}
\usepackage{euscript}
\usepackage{fancybox}
\usepackage{graphicx}
\usepackage{hyperref}
\usepackage{lipsum}
\usepackage{pgf}
\usepackage{psfrag}
\usepackage{pifont}
\usepackage{mathrsfs}
\usepackage{multirow}
\usepackage{stfloats}
\usepackage{setspace}
\usepackage{stfloats}
\usepackage{tabularx}
\usepackage{textcomp}

\newtheorem{Theorem}{Theorem}

\newtheorem{definition}{Definition}
\newtheorem{Remark}{Remark}
\SetKwFunction{Dothat}{Do that}

\newcolumntype{C}[1]{>{\centering\arraybackslash}p{#1}}

    \def\Complex{{\rm\rule[.23ex]{.03em}{1.1ex}\kern-.3em{C}}}

    \newcommand{\be}{\begin{equation}} \newcommand{\ee}{\end{equation}}
    \newcommand{\bea}{\begin{eqnarray}} \newcommand{\eea}{\end{eqnarray}}
    \newcommand{\benum}{\begin{enumerate}} \newcommand{\eenum}{\end{enumerate}}

    \newcommand{\qw}{{\bf w}}

    \newcommand*{\argmin}{\operatornamewithlimits{argmin}\limits}
    \newcommand*{\argmax}{\operatornamewithlimits{argmax}\limits}

\begin{document}

\title{Intelligent Charging Management of Electric Vehicles Considering Dynamic User Behavior and Renewable Energy: A Stochastic Game Approach}

\author{ \IEEEauthorblockN{Hwei-Ming Chung, \IEEEmembership{Student Member, IEEE}, Sabita Maharjan, \IEEEmembership{Senior Member, IEEE}, Yan Zhang, \IEEEmembership{Fellow, IEEE}, and Frank Eliassen, \IEEEmembership{Member, IEEE} }

\thanks{This work was supported by Norwegian Research Council under Grants 275106 (LUCS project), 287412 (PACE project), and 267967 (SmartNEM project).
}

\thanks{H.-M. Chung and F. Eliassen are with the Department of Informatics, University of Oslo, Oslo 0373, Norway, (e-mail: hweiminc@ifi.uio.no, frank@ifi.uio.no).}

\thanks{S. Maharjan and Y. Zhang are with Department of Informatics, University of Oslo, Oslo 0373, Norway; and Simula Metropolitan Center for Digital Engineering, Oslo 0167, Norway. (e-mail: yanzhang@ieee.org, sabita@ifi.uio.no).
}

}

\maketitle
\begin{abstract}
Uncoordinated charging of a rapidly growing number of electric vehicles (EVs) and the uncertainty associated with renewable energy resources may constitute a critical issue for the electric mobility (E-Mobility) in the transportation system especially during peak hours.
To overcome this dire scenario, we introduce a stochastic game to study the complex interactions between the power grid and charging stations.
In this context, existing studies have not taken into account the dynamics of customers' preference on charging parameters.
In reality, however, the choice of the charging parameters may vary over time, as the customers may change their charging preferences. 
We model this behavior of customers with another stochastic game.
Moreover, we define a quality of service (QoS) index to reflect how the charging process influences customers' choices on charging parameters.
We also develop an online algorithm to reach the Nash equilibria for both stochastic games.
Then, we utilize real data from the California Independent System Operator (CAISO) to evaluate the performance of our proposed algorithm.
The results reveal that the electricity cost with the proposed method can result in a saving of about $20\%$ compared to the benchmark method, while also yielding a higher QoS in terms of charging and waiting time.
Our results can be employed as guidelines for charging service providers to make efficient decisions under uncertainty relative to power generation of renewable energy.

\end{abstract}

\section*{Nomenclature} 
\addcontentsline{toc}{section}{Nomenclature}
\begin{IEEEdescription}[\IEEEusemathlabelsep\IEEEsetlabelwidth{$V_1,V_2$}]
\item[A. Sets and Indices]
\item[$i$] EV index.
\item[$m$] Charging station index.
\item[$t$] Time index.
\item[${\cal T}$]  Time horizon.
\item[$M$]  Number of charging stations.
\item[${\cal H}_{m, t}$] Set of EVs in the $m$-th charging station.
\\

\item[B. Variables]
\item[${\tt P}_{i, t}$]  Charging rate of EV $i$.
\item[$x_{m, t}$]  Power purchased from the grid for charging station $m$ at time $t$.
\item[$\omega_{i, j}$] Probability of EV $i$ choosing the $j$-th deadline.
\\

\item[C. Parameters]
\item[${\tt SOC}_{i, t}$]  SOC of EV $i$ at time $t$.
\item[${\tt SOC}_{i}^{\rm fin}$]  Target SOC of EV $i$.
\item[${\tt E}_{i}^{\tt cap}$]  Battery capacity of EV $i$.
\item[${\tt E}_{i}^{\tt de}$]  Demand of EV $i$.
\item[${\tt E}_{i}^{\tt fin}$]  Energy level of the battery when EV $i$ leaves. 
\item[${\tt P}_{i}^{\tt max}$] Maximum charging rate of EV $i$.
\item[$\epsilon_{c}$] Charging efficiency of the EV.
\item[$v_{i}$ ($a_{i}$)] Departure (arrival) time of EV $i$.
\item[$f_{i}$]  Estimated finishing time of the charging of EV $i$.
\item[$r_{m, t}$]  Renewable power generation of the $m$-th charging station at time $t$.
\item[${\tt P}_{t}^{{\tt base}}$] Based load of the power grid at time $t$.
\item[${\tt P}^{{\tt peak}}$] Peak constraint of the power grid.
\item[$k_{t}$] Electricity price at time $t$.
\item[$Q_{m, t}$] Charging queue of charging station $m$ at time $t$.
\item[$Z_{m, t}$] Request penalty queue of charging station $m$ at time $t$.
\item[$B_{i, t}$] Waiting time queue of EV $i$ at time $t$.
\\

\item[D. Operators]
\item[$|\cdot|$]   Cardinality of set.
\item[$\mathbf{A}^{T}$]   Transpose of matrix $\mathbf{A}$.

\item[Other notations are defined in the text.]  

\end{IEEEdescription}

\section{Introduction} \label{sec:introduction}

Electrification of transportation (E-Mobility) plays an important role in achieving the goals of the Paris agreement enforced in $2016$, which requires $1170$ giga tons of $\mbox{CO}_{2}$ from cumulative emissions to be cut down during the $2015-2100$ period \cite{2017-EV-outlook}.
Under this trend, the deployment of EVs is regarded as an important initiative to reduce $\mbox{CO}_{2}$ emission from conventional vehicles.
Various policies have been made in many countries in order to reach the vehicle market share of EVs up to $30\%$ by $2030$ \cite{2018-EV-outlook}.

One of the key issues of deploying EVs is the driving distance.
However, with a large and rapidly growing number of EVs, charging service providers in the transportation system should be aware of the negative impact of uncoordinated charging to the efficient and reliable operation of the power grid.
Uncoordinated large-scale EV charging can lead to problems, such as increase in peak load and power quality reduction \cite{2015-ev-impact}.
To cope with this situation, researchers have proposed various coordinated charging mechanisms (e.g., \cite{2012-He-EV-selection,2018-shi-EV-selec,2018-david-EV-manage,2019-Seyedyazdi-EV-manage,2019-zhang-EV-manage,2012-wen-EV-selection,2017-mal-EV-selection,2018-chung-EV-selection}).
These algorithms require various charging parameters, including charging demand, arrival time, and departure time, as inputs.

The EV charging management problem was studied in \cite{2012-He-EV-selection,2018-shi-EV-selec,2018-david-EV-manage,2019-Seyedyazdi-EV-manage,2019-zhang-EV-manage}, where the interaction between power grid and charging stations was considered.  
The authors in \cite{2012-He-EV-selection} presented a mathematical model in order to minimize electricity cost for charging EVs.
Based on the framework in \cite{2012-He-EV-selection}, the authors in \cite{2018-shi-EV-selec} applied model predictive control to jointly address optimal power flow problem for managing charging tasks.
Then, the competition between charging stations may lead to a less efficient charging management solution.
Therefore, irrational behavior of the charging stations was taken into account in \cite{2018-david-EV-manage}.
The authors in \cite{2019-Seyedyazdi-EV-manage} further studied the dynamic behaviors of the charging stations and EV owners.
The method of optimally utilizing the charging capacity was proposed in \cite{2019-zhang-EV-manage}.
On the other hand, studies such as \cite{2012-wen-EV-selection,2017-mal-EV-selection} schedule the charging of EVs without considering the power generation from the grid.
The authors in \cite{2012-wen-EV-selection,2017-mal-EV-selection} proposed a parameter called user convenience, which is calculated from charging states of EVs, and algorithms were proposed to maximize user convenience.
The authors in \cite{2018-chung-EV-selection} further applied user convenience to the charging management problem.

Integration of renewable energy resources into the transportation system is another initiative in order to further reduce $\mbox{CO}_{2}$ emission.
However, dealing with inherent uncertainty associated with power generation from renewable energy resources, such as solar panels and wind turbines, will be a major issue for the charging management problem.
Such uncertainties were addressed in \cite{2015-zhao-stochastic-infocom, 2016-wayes-EV-manage,2018-luo-stochastic,2017-huang-stochastic, 2018-huang-stochastic,2018-yang-stochastic,2018-hiroshi-stochastic,2018-zhou-EV-lyapunov}.
In \cite{2015-zhao-stochastic-infocom}, the authors proposed an index, called competitive ratio, to design an algorithm that can capture system uncertainty; however, the competitive ratio needs the optimal offline solution.
The authors in \cite{2016-wayes-EV-manage} classified EVs with different priorities to manage the uncertainty associated with renewable energy generation.
The multi-objective optimization was formulated to solve the EV charging management problem with renewable energy in \cite{2018-luo-stochastic}.
The authors in \cite{2017-huang-stochastic,2018-huang-stochastic,2018-yang-stochastic} modeled the EV charging management problem as Markov decision process (MDP).
The uncertainty associated with renewable power was addressed in both power flow operation and charging management problem by the authors in \cite{2018-hiroshi-stochastic}.
The charging management problem in a charging station was formulated as a stochastic optimization problem in \cite{2018-zhou-EV-lyapunov}.

For the studies \cite{2017-huang-stochastic,2018-huang-stochastic,2018-yang-stochastic,2018-hiroshi-stochastic}, they showed that EV charging management problem with the uncertainty of renewable energy generation can be solved with MDP and obtained promising results in terms of electricity cost.
However, for utilizing MDP, we need to model transition probabilities of the system states (e.g., renewable energy generation or charging demand), which are not readily available in practice.
Thus, we rely on predicted information.
Unfortunately, such predictions have often proven to be quite inaccurate \cite{2018-gigoni-forecast-error}.
For example, the day-ahead prediction error in some days for renewable energy generation may result in a large error due to the unpredictable weather \cite{2018-gigoni-forecast-error}.
Therefore, designing an online algorithm for EV charging management incorporating renewable energy resources remains a big challenge and an important problem to solve.
Also, different from \cite{2015-zhao-stochastic-infocom,2016-wayes-EV-manage,2018-zhou-EV-lyapunov}, we study EV charging management problem among several charging stations that requires cooperation between charging stations.
In addition, the dynamic behavior of EV owners or charging stations during charging should be also considered.
In \cite{2018-david-EV-manage}, irrational behaviors of the charging stations were modeled by using prospect theory.
Prospect theory was also used to represent customers' irrationally varying willingness of participating in a demand response program \cite{2016-wang-DR-user}. 
Furthermore, the driving behaviors were incorporated in the driving range estimation in \cite{2014-drive-range}.
The preferences of EV owners on choosing charging stations and the preferences of charging stations on choosing EV owners were jointly considered in \cite{2019-Seyedyazdi-EV-manage}.
However, the dynamic behavior of EV owners with regard to the charging parameters also plays an important role in the charging management problem.
Specifically, algorithms designed in \cite{2012-He-EV-selection,2018-shi-EV-selec,2012-wen-EV-selection,2017-mal-EV-selection,2018-chung-EV-selection} are based on the parameters given by EV owners.
However, EV owners may change preference on the charging parameters during charging and they will not explicitly report this dynamic behavior to the charging stations.
Moreover, the dynamic preference on the charging parameters may create conditions detrimental to reliability and resiliency of the charging algorithms.
Therefore, when designing the algorithm for solving the charging management problem, it should be taken into account that EV owners may change the charging parameters while the EVs are in the charging stations.
Jointly considering and modeling the dynamics of customers' preference on charging parameters and uncertainty associated with renewable energy resources, therefore demand novel solutions.

To overcome the aforementioned challenges, in this paper, we employ a stochastic game, a generalization of MDP, such that we rely to a lesser extent on the forecasted information.
The forecasted information is only used to calculate the charging status of EVs for EV owners.
Also, EV owners may dynamically change their preference on charging parameters.
To capture and characterize this dynamic behavior of EV owners, we introduce a QoS index.
More specifically, EV owners attempt to obtain a better QoS during the charging process by changing their preferences on charging parameters.
We then model varying choices of EV owners on charging parameters under the uncertainty associated with renewable power generation also as a stochastic game.
In order to obtain the Nash equilibria for both stochastic games, we design an online algorithm.
We evaluate the performance of our proposed solution using real-world data from CAISO \cite{Load_data}, that provides all the historical data of the power system operation in California.
We show that the proposed algorithm yields satisfactory results even if the prediction error is large.

The main contributions of this paper are threefold:
\begin{itemize}

\item In the paper, we address the EV charging management problem incorporating the dynamic behavior of EV owners in terms of changing the charging parameters.
We model the problem in the form of two stochastic games where the amount of renewable power generated is the stochastic variable.

\item We design an online algorithm to approach the Nash equilibria for both stochastic games that leads to a suboptimal solution.
We further quantify the proximity of the obtained solution to the optimal solution.

\item We utilize real-world data to illustrate the performance of the proposed algorithm in terms of total electricity cost and QoS.
The results show a reduction of about $20\%$ in electricity cost and a gain of about $12\%$ in terms of QoS.

\end{itemize}


\section{System Model}\label{sec:system_model}

\begin{figure}
\begin{center}
\resizebox{3.2in}{!}{%
\includegraphics*{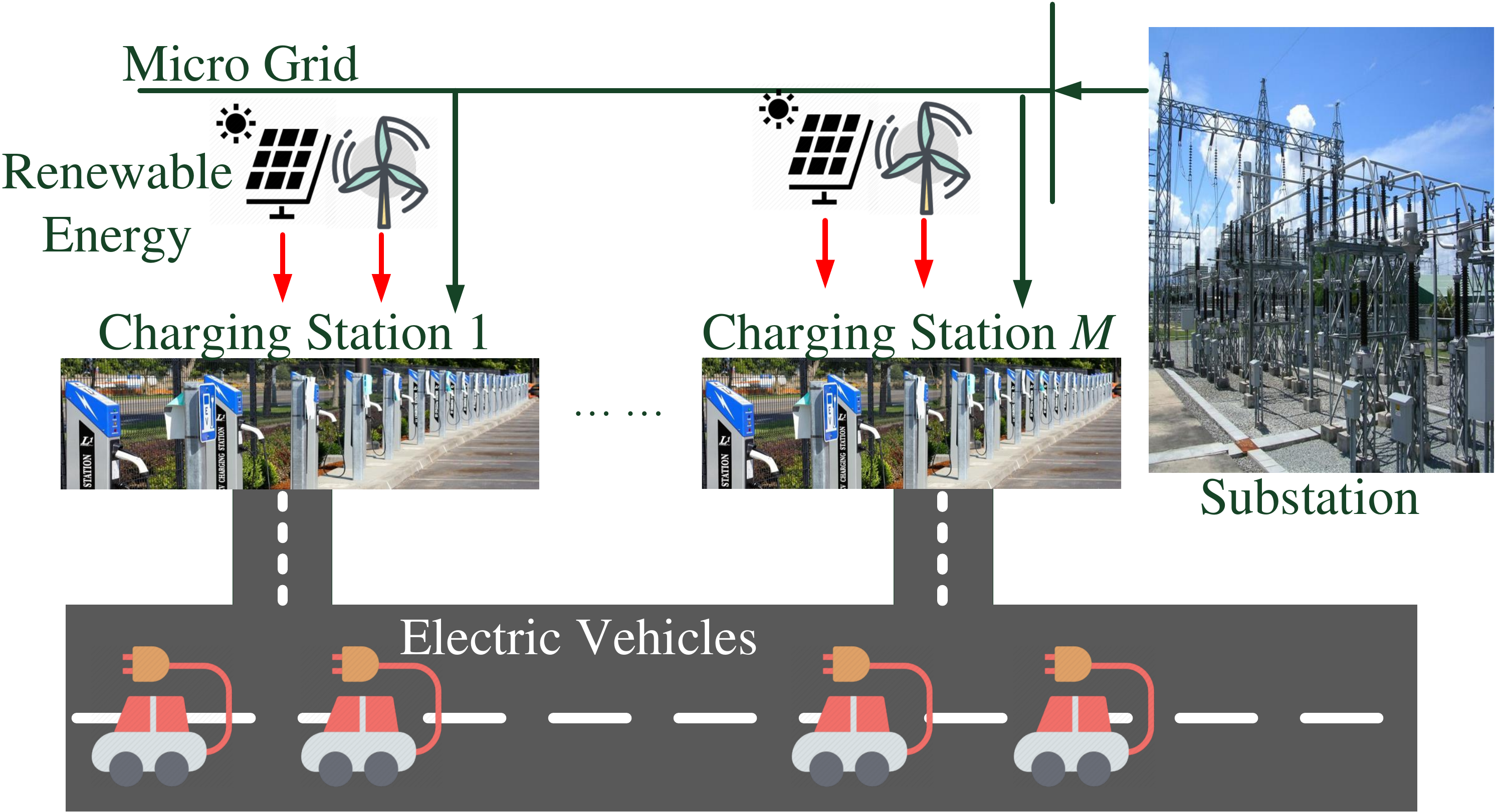} }%
\caption{System model used in this paper.} 
\label{fig:System_Model}
\end{center}
\end{figure}

In this paper, we consider a transportation system consisting of $M$ charging stations each having renewable energy resources (e.g., solar or wind power) as illustrated in Fig. \ref{fig:System_Model}.
All charging stations are connected to power grid and belong to a single utility or single private enterprise.
Then, the owner of the charging stations collects the charging parameters from EV owners and then schedule the charging tasks.
We consider a total of $N$ EVs in the system.
The charging management problem is studied in a time horizon ${\cal T}$ with equal length time slots $t=1, 2,\ldots,{\cal T}$.
In time slot $t$, ${\cal H}_{m, t}$ denotes the set of EV indices in the $m$-th charging station.
Moreover, we assign EV owners the same indices as their EVs.
Then, ${\cal H}_{m, t}$ is used to represent the cardinality of the set ${\cal H}_{m, t}$.
We use ${\cal H}_{t}= \bigcup^{M}_{m=1} {\cal H}_{m, t}$ to denote the set of EV indices in all charging stations.

For EV $i$, the state of charge (SOC) at time slot $t$ is denoted as ${\tt SOC}_{i, t}$, where $0 \leq {\tt SOC}_{i, t} \leq 1$.
Then, the capacity of EV $i$ is denoted as ${\tt E}_{i}^{\tt cap}$.
The target SOC, denoted as ${\tt SOC}_{i}^{\rm fin}$, indicates the expected SOC at the departure time of EV $i$. 
We further employ ${\tt E}_{i}^{\tt de}$ and ${\tt E}_{i}^{\tt fin}$ to denote demand of EV owner $i$ and the final battery energy level at departure time, respectively.
Then, $a_i$ denotes the arrival time of EV $i$, and $v_{i}$ is the corresponding deadline (departure time).
With these notations, we can write ${\tt E}_{i}^{\tt de} = {\tt E}_{i}^{\tt fin} - {\tt SOC}_{i, a_i}  {\tt E}_{i}^{\tt cap}$.
When EV $i$ arrives in the charging station, EV owner submits the charging parameters (i.e., $a_{i}$, $v_{i}$, and ${\tt SOC}_{i}^{\rm fin}$) to the charging station.

In every time slot, EV $i$ is assigned a charging rate, denoted as ${\tt P}_{i, t}$, which is restricted by the maximum charging rate, ${\tt P}_{i}^{\tt max}$, based on the physical constraints of batteries.
Other than the charging rate, the charging station can broadcast the estimated finishing time slot of the charging task, denoted as $f_{i}$, to EV owner $i$.
Also, $f_{i}$ is always earlier than $v_{i}$.

In the $m$-th charging station, the generation of renewable energy at time $t$ is represented as $r_{m, t}$, and its capacity is denoted as $r_{m}^{{\tt cap}}$.
However, the available renewable energy may not be enough to meet the demand of all EVs.
To overcome this situation, the $m$-th charging station has to purchase an amount of energy, $x_{m, t}$, from the grid.
The base load and the peak constraint of the power grid at time $t$ are ${\tt P}_{t}^{{\tt base}}$ and ${\tt P}^{{\tt peak}}$, respectively.
Thus, we have 
\begin{equation}
{\tt P}_{t}^{{\tt base}} + \sum_{m=1}^{M} x_{m, t} \leq {\tt P}^{{\tt peak}}.
\end{equation}

Throughout this paper, we consider the following:
\begin{enumerate}
\item the charging station only has knowledge about the base load and the generation of renewable energy at the current time slot.
The base load information and the generation of renewable energy in the future can only be forecasted; and

\item EVs do not discharge power to the grid (i.e., ${\tt P}_{i, t} \leq 0$), since discharging infrastructure for the charging stations is still expensive and discharging may accelerate battery degradation rate; and

\item real-time pricing model is adopted.
Therefore, charging stations can only know the current electricity price.
\end{enumerate}

To measure how the EV owners act accordingly to the charging processes, a QoS index is required.
Some charging algorithms may let every EV spend long time on charging, such that all EV owners wait long time for SOC reaching ${\tt SOC}_{i}^{fin}$.
On the other hand, some algorithms select EVs with priority, and therefore the waiting time of the EVs may vary considerably.
Hence, QoS for EV $i$ can be defined as
\begin{equation} \label{eq:QoS}
\mbox{QoS}_{i} =  \ln (1 + (f_{i} - a_{i})) + \ln \left(1 + \left( v_{i} - f_{i} \right) \right). 
\end{equation}
The first term indicates the charging time, which is the duration for the SOC to reach ${\tt SOC}_{i}^{\rm fin}$.
The second term represents the duration after SOC reaching ${\tt SOC}_{i}^{\rm fin}$ but before the deadline.
Normally, the total time for EVs staying in the charging stations can be regarded as the summation of charging and waiting time.
Then, the power for charging EVs is determined by the charging stations, and therefore the charging stations can be regarded as the entities that allocate the charging and waiting time for EVs.
The EV owners can evaluate the charging process with (\ref{eq:QoS}), and then adjust the charging parameters to influence the power allocation of the charging stations such that the QoS can be improved.
The definition of QoS is motivated by the objective function of the fair rate allocation problem \cite{1998-rate-allocation}, where QoS in (\ref{eq:QoS}) represents the fairness of the charging process.


\section{Stochastic Game Formulation}\label{sec:game_model}

As shown in Fig. \ref{fig:System_Model}, there exists a two-level interaction between three participants in E-Mobility, namely substation, charging stations, and EV owners.
The first level is between the substation and the charging stations, and the second level is between the charging station and the EV owners.
In the first level, the charging stations decide the amount of power to draw from the power grid.
Then, based on the available power, the charging station can broadcast the estimated finishing charging time to EV owners.
After receiving $f_{i}$, EV owners may want to change the charging parameters to improve QoS.
Because the uncertainty of renewable energy is a stochastic variable, which can influence the decisions of the participants, we use stochastic games to model interactions of both levels.
Specifically, we game the first level as a cost minimization game and game the second level as a deadline selection game.

\subsection{Cost Minimization Game} \label{subsec:charging_game}

At the beginning of every time slot, the charging stations have to decide the amount of power to draw from the power grid.
For achieving better E-Mobility, charging stations require more power from the power grid than before; however, some uncertainties may exist so that it is difficult for charging stations to make decisions.
Renewable power generation in the future and future electricity prices are uncertain.
Also, substations may increase the electricity price in the upcoming time slot if charging stations consume too much power at current time slot.
To minimize the total electricity cost in a day, charging stations should interact with the substations.
This interaction can be modeled as a cooperative stochastic game, which can be defined by a tuple  $\Xi_{1} = \{ k_{t}, \mathbf{r}, \mathbf{x}, \{1, 2, \cdots, M \}, u_{m}\}$.
The main components of the game include:
\begin{itemize}
\item $k_{t}$ represents real-time electricity price at time $t$; 
\item $\mathbf{r} := \{r_{1, t}, \cdots, r_{M, t}\}$ represents states of the renewable power generation in all charging stations; 
\item $\mathbf{x} := \{x_{1, t}, \cdots, x_{M, t}\}$ is the action of the charging stations;
\item $\{1, 2, \cdots, M \}$ is the set of the players (i.e., charging stations in the power grid); and
\item $u_{m}$ is the expected payoff function of the $m$-th charging station corresponding to the selected action. 
\end{itemize}
The $m$-th charging station observes the amount of renewable energy generated, $r_{m, t}$, and the real-time electricity price, $k_{t}$, at each time slot.
Then, the charging stations have to decide how much power to purchase from the grid (i.e., $\mathbf{x}$).
After selecting the amount of energy to purchase from the grid, the charging stations receive a payoff, $u_{m}$.
After selecting $\mathbf{x}$, the charging stations receive a payoff, $u_{m}$.
The payoff function for the charging stations is the expected electricity cost in a day, which can then be expressed as 
\begin{equation}
u_{m}(x_{m, t}, \mathbf{x}_{-m, t}) = \mathop{\mathbb{E}} \left[\sum_{t=1}^{ {\cal T} } k_{t}  \sum_{m=1}^{M} x_{m, t}  \right],
\end{equation}
where $\mathbf{x}_{-m, t}$ is the power drawn by the charging stations except the $m$-th charging station.

In the cost minimization game, the charging stations aim to minimize the expected electricity cost in a day; we can search for the Nash equilibrium (NE) to reach this goal.
NE for the cost minimization game is defined as following:

\begin{definition}[NE for the cost minimization game] \label{def:NE_def}
In game theory, NE is a state of the game that no player can benefit by unilaterally changing strategies, while the other players keep their strategies unchanged.
The NE of the cost minimization game $\Xi_{1}$, denoted as $\mathbf{x}^{*}:= \{ x_{1, t}^{*} \cdots x_{M, t}^{*} \}$, such that for all charging stations, it satisfies the following inequality:
\begin{equation}
u_{m}( x^{*}_{m, t}, \mathbf{x}^{*}_{-m, t}) \leq u_{m}( x_{m, t}, \mathbf{x}^{*}_{-m, t}),
\end{equation}
where it implies that the electricity cost will be higher if the $m$-th charging station does not choose $x^{*}_{m, t}$.
\end{definition}

The NE of the cost minimization game can be obtained by solving the following optimization problem.
\begin{subequations} \label{eq:NE_power_decision_game} 
\begin{align}
\min_{\mathbf{x}} &  ~~u_{m}(x_{m, t}, \mathbf{x}_{-m, t})  & \\
\mbox{subject to}    & ~~ \mathop{\mathbb{E}} \left[ \sum_{t} \left( x_{m, t} + r_{m, t} \right) \right] \geq \sum_{i \in {\cal H}_{m, t} } {\tt E}_{i }^{\rm de}, & \forall m  \label{eq:meet_de_expect} \\
                     & ~~{\tt P}_{t}^{{\tt base}} + \sum_{m=1}^{M}  x_{m, t}  \leq {\tt P}^{{\tt peak}}, &  \forall t. \label{eq:lower_peak_expect}          
\end{align} 
\end{subequations}
Constraint (\ref{eq:meet_de_expect}) indicates that the amount of power purchased from the power grid and the generated renewable power should fulfill the total charging demand.
Constraint (\ref{eq:lower_peak_expect}) represents the capacity constraint, i.e., the total load in a power grid should be bounded by ${\tt P}^{{\tt peak}}$.
The uncertainty associated with the renewable power level is regarded as a random variable and the total electricity cost should be calculated under this uncertainty.
Therefore, (\ref{eq:NE_power_decision_game}) is a stochastic optimization problem.

\subsection{Deadline Selection Game} \label{subsec:deadsec_game}

After deciding $x_{m, t}$, the charging station allocates power to EVs and then estimates $f_{i}$ for EV owners.
The generation of the renewable power in the upcoming time slots may influence the outcome of $f_{i}$, and therefore $f_{i}$ can be regarded as a stochastic term to EV owners.

After receiving $f_{i}$ from the charging station, the EV owners can calculate the waiting time.
According to (\ref{eq:QoS}), the QoS may be very low for some EV owners.
Then, from the customers' perspective, the EV owners may change their preferences on the charging parameters, such that the power allocation to EVs can be influenced.
By doing so, the QoS can be improved.
For the charging parameters, $a_{i}$ is fixed when the EV arrives, and therefore it cannot be modified.
Thus, EV owners can change $v_{i}$ and ${\tt SOC}_{i}^{\rm fin}$ according to the charging process.
In our model, EV owners can change $v_{i}$ in order to influence QoS.
This kind of stochastic interactions has not been incorporated in previous studies.
Similar to the cost minimization game, the interaction between EV owners and the charging station can be formulated as a non-cooperative stochastic game called deadline selection game.
The game can also be represented as a tuple, $\Xi_{2} = \{ \mathbf{f}, \mathbf{D}, {\cal H}_{m, t}, u_{i} \}$, and the main components in this stochastic game include:
\begin{itemize}
\item $\mathbf{f} := \left \{f_{1}, \cdots, f_{ |{\cal H}_{m, t}|} \right \}$ denotes the estimated finishing time of charging broadcasted by the charging station;
\item $\mathbf{D} := \{ {\bf d}_{1}^{T}, \cdots, {\bf d}_{ |{\cal H}_{m, t}|}^{T} \}$ is the action space for the EV owners;
\item ${\cal H}_{m, t}$ is the set of the players (i.e., EV owners with EVs in the $m$-th charging station); and
\item $u_{i}$ is the payoff function to the EV owner $i$.
\end{itemize}
The charging station will broadcast the estimated finishing charging time, $f_{i}$, to the EV owner $i$ after allocating the power to EVs, and $f_{i}$ is defined as
\begin{equation}
 f_{i} = \left\{ ~ t ~ \bigg\vert ~~{\tt SOC}_{i, t} = \frac{{\tt E}_{i}^{\tt fin}}{{\tt E}_{i}^{\tt cap}} \right\}.
\end{equation}
After receiving $f_{i}$ from the charging station, the EV owner $i$ can take $k$ different actions in the action space denoted as ${\bf d}_{i}^{T} :=  \{  d_{i, 1} \cdots d_{i, k}\}$, where the elements in the action space revolve around $v_{i}$.
Then, let $\mathbf{W} := \{ \pmb{\omega}^{T}_{1}, \cdots, \pmb{\omega}^{T}_{ |{\cal H}_{m, t}|} \}$ denote the vector of the mixed strategies of all EV owners, where, for each EV owner $i \in {\cal H}_{m, t}$, we have $\pmb{\omega}_{i}^{T} := \{ \omega_{i, 1} \cdots \omega_{i, k}\}$ and $\omega_{i, j}$ is the probability corresponding to the choice of the $j$-th pure strategy in ${\bf d}_{i}^{T}$.
Based on the response, the EV owners will receive a payoff, which is the risk function denoted as
\begin{equation} \label{eq:EV_payoff}
\!\!\!\! u_{i} \left( \pmb{\omega}_{i}^{T}, \pmb{\omega}_{-i}^{T} \right) \!=\! \sum_{i \in {\cal H}_{m, t}} \mathop{\mathbb{E}} \!\!\left[ \sum_{j=1}^{k} w_{i, j} ( {\tt SOC}_{i}^{\rm fin} - {\tt SOC}_{i, t}) e^{ f_{i} - d_{i, j} } \right],
\end{equation}
where $\pmb{\omega}_{-i}^{T}$ denotes the mixed strategy of the EV owners except the $i$-th EV owner. 
Here, the risk function indicates the risk level of fulfilling the demand under the chosen deadline for the EV owner.
The goal of EV owner $i$ is to choose a mixed strategy so as to minimize its payoff function as given in (\ref{eq:EV_payoff}).

The NE for the deadline selection game can be obtained by solving the following optimization problem
\begin{subequations} \label{eq:NE_deadline_selection_game} 
\begin{align}
\min_{ \mathbf{W} } &  ~~ u_{i} \left( \pmb{\omega}_{i}^{T}, \pmb{\omega}_{-i}^{T} \right)  & \\
\mbox{subject to}    & \quad  \sum_{j=1}^{k} w_{i, j} = 1, \label{eq:NE_weight_sum_per_EV}  & \forall i \in {\cal H}_{m, t} \\
                     &  \quad w_{i, j}  \geq 0. \label{eq:NE_weight_cons}     &  \forall i \in {\cal H}_{m, t}, j 
\end{align} 
\end{subequations}
The objective function of (\ref{eq:NE_deadline_selection_game}) aims to minimize the risk function for the EV owners.
Constraint (\ref{eq:NE_weight_sum_per_EV}) indicates that the sum of all probabilities over the action space should be $1$.
Constraint (\ref{eq:NE_weight_cons}) ensures that $w_{i, j}$ is non-negative.
The uncertainty about renewable power generation may influence the estimation on $f_{i}$, and therefore it is also a stochastic variable
Thus, (\ref{eq:NE_deadline_selection_game}) is also a stochastic optimization problem.

\section{Equilibria for the Games}\label{sec:problem_trans}
From (\ref{eq:NE_power_decision_game}) and (\ref{eq:NE_deadline_selection_game}), it is clear that obtaining the Nash equilibria is challenging due to the expectation term.
Also, there are many uncertainties such as the electricity prices, states of the renewable power, and the demand of upcoming EVs.
To avoid the need of predicting such parameters, we will design an online algorithm, for which we need to convexify the original problems.
We use virtual queues \cite{2010-michael-lyapunov} to model the charging process and then design an online algorithm based on the virtual queue.
The steps are explained in detail in this section.

\subsection{Cost Minimization Game Convexification} \label{subsec:power_decision_queue}

While searching the NE, the charging station has to meet the demand from EV owners as stated in (\ref{eq:meet_de_expect}).
Therefore, for the $m$-th charging station, we can construct a virtual queue to address the current charging tasks.
The charging virtual queue can be expressed as
\begin{equation} \label{eq:Q-virtual-queue}
Q_{m, t+1} = \max\{Q_{m, t} - \epsilon_{c} Y_{m, t}, 0 \} + \lambda_{m, t},
\end{equation}
where $Y_{m, t} = x_{m, t} + r_{m, t}$ is the total available energy for charging EVs and $\lambda_{m, t}$ is the demands of the EVs arriving at time $t$.
Quantity $\epsilon_{c}$ is the charging efficiency.

Note that (\ref{eq:Q-virtual-queue}) does not include the terms accounting for the deadline of charging tasks.
Thus, we create another virtual queue for the $m$-th charging station called request penalty queue to address this issue, i.e., 
\begin{equation}\label{eq:Z-virtual-queue}
Z_{m, t+1} = \max\{Z_{m, t} +  \eta_{m} Q_{m, t} - \epsilon_{c}Y_{m, t}, 0\},
\end{equation}
where $\eta_{m}$ is the parameter that can adjust the growth rate.

The Lyapunov function associated with the virtual queues constructed above will be
\begin{equation} \label{eq:first_lyapunov_func}
L_{t}^{\tt charg} = \frac{1}{2} \left( \sum_{i=1}^{M}  Q_{m, t}^{2}  + \sum_{i=1}^{M}  Z_{m, t}^{2}  \right).
\end{equation}
Then, we define a one-step conditional Lyapunov drift as
\begin{equation} \label{eq:charging_drift}
\Delta_{t}^{\tt charg} = \mathop{\mathbb{E}} \left[ L_{t+1}^{\tt charg} - L_{t}^{\tt charg} | Q_{m, t}, Y_{m, t}, \lambda_{m, t} \right],
\end{equation}
where it represents the expected change in the Lyapunov function from $t$ to $t+1$.
In (\ref{eq:charging_drift}), it is difficult to capture the probability density function of the renewable energy.
Therefore, instead of minimizing (\ref{eq:charging_drift}) directly, we can minimize its upper bound.
The individual upper bound of the virtual queues will be
\begin{subequations}
\begin{align}
&  \!\!\!\! Q_{m, t+ 1}^{ 2 }  \!-\! Q_{m, t}^{2}  \leq  \epsilon_{c}^{2} Y_{m, t}^{2} + 2 Q_{m, t} \left( \lambda_{m, t}^{\tt max} -  \epsilon_{c}Y_{m, t} \right),  \label{eq:upper_Q} \\
&  \!\!\!\! Z_{m, t+1}^{2}  \!-\!  Z_{m, t}^{2}  \!\leq \! 2 Z_{m, t} \left(\eta_{m} Q_{m, t} \!-\! \epsilon_{c} Y_{m, t} \right)  +  ( \eta_{m} Q_{m, t} \!-\! \epsilon_{c} Y_{m, t} )^{2}.\!\!\! \label{eq:upper_Z}
\end{align}
\end{subequations}
Based on (\ref{eq:upper_Q}) and (\ref{eq:upper_Z}), Lyapunov drift can be bounded by
\begin{align}
\Delta_{t}^{\tt charg} & \leq \sum_{m=1}^{M} \mathop{\mathbb{E}} \! \left[ \frac{1}{2} \left( \epsilon_{c}^{2} Y_{m, t}^{2}  +  (\eta_{m} Q_{m, t} - \epsilon_{c} Y_{m, t} )^{2}  \right)\right]    \notag \\
&  + \sum_{m=1}^{M} \mathop{\mathbb{E}} \!\left[ Q_{m, t}\! \left( \lambda_{m, t}^{\tt max^{2}} \!-\! \epsilon_{c} Y_{m, t} \right) \!+\! Z_{m, t}  \left(\eta_{m} Q_{m, t} \!-\! \epsilon_{c} Y_{m, t} \right)  \right]  \notag \\
& =  \mathop{\mathbb{E}} \left[J_{1} (\mathbf{x} )| Q_{m, t}, Y_{m, t}, \lambda_{m, t} \right] \label{eq:Q_final_eq}
\end{align}
where we regard $\lambda_{m, t}^{\tt max^{2}}$ as a constant in (\ref{eq:Q_final_eq}).
The detailed derivation is provided in Appendix \ref{subsec:first_queue}.

Minimizing $J_{1} (\mathbf{x})$ ensures that $\mathop{\mathbb{E}} \left[J_{1}(\mathbf{x}) | Q_{m, t}, Y_{m, t}, \lambda_{m, t} \right]$ is minimized by opportunistically minimizing an expectation \cite{2010-michael-book}.
However, minimizing the Lyapunov drift upper bound, $\Delta_{t}^{\tt charg}$, can push the queue to a lower value.
Therefore, we further add another penalty term to $\Delta_{t}^{\tt charg}$, such that the stability of the queues can be considered.
Then, this is called \textit{drift-plus-penalty} algorithm.
The \textit{drift-plus-penalty} algorithm states that we should greedily minimize the upper bound of the \textit{drift-plus-penalty}.
Hence, we formulate an optimization problem based on the \textit{drift-plus-penalty} algorithm as 
\begin{subequations} \label{eq:P1}
\begin{align}
\mathcal{P}_1:  \min_{ \mathbf{x} } &  \quad J_{1} (\mathbf{x} )  + k_{t}  \left( \sum_{m=1}^{M} V_{m}^{\tt charg}  x_{m, t}  \right) & \label{eq:P1_obj} \\
\mbox{subject to}    & ~~Y_{m, t} \leq \min \left\{ \sum_{i \in {\cal H}_{m, t}} {\tt P}_{i}^{\tt max}, Q_{m, t} \right\}, & \forall m \label{eq:not_over_charging_limit} \\
                     & ~~{\tt P}_{t}^{{\tt base}} + \sum_{m=1}^{M} x_{m, t} \leq {\tt P}^{{\tt peak}}. \label{eq:not_over_peak}      &    
\end{align} 
\end{subequations}
Here, the penalty term is the total electricity cost at time $t$, and $V_{m}^{\tt charg}$ is the weight with respect to the electricity cost for the $m$-th charging station.
Thus, (\ref{eq:P1_obj}) is a trade-off between minimizing the queue-length drift and the penalty function.
Some charging stations may not receive enough power in the current time slot, and therefore we put the remaining demand into the penalty queue.
If the electricity price is the same in the next time slot, the charging station with huge value in the penalty queue should receive more power in the next time slot. 
This is because we minimize the expected difference of the summation of two virtual queues as in (\ref{eq:charging_drift}).
The optimization problem $\mathcal{P}_1$ is a quadratic optimization problem, which can be solved using interior-point method \cite{boyd-cvx-book}. 
Constraint (\ref{eq:not_over_charging_limit}) ensures that $Y_{m, t}$ cannot exceed the amount that all EVs are charged with ${\tt P}_{i}^{\tt max}$ or the current queue value.
Eq.\,(\ref{eq:not_over_peak}) keeps the base load with the power drawn from the grid less than the peak constraint, ${\tt P}^{{\tt peak}}$.

\subsection{Deadline Selection Game Convexification} \label{subsec:dead_selec_queue}

In the deadline selection game, the EV owners attempt to change the deadline, while the risk level is considered.
Then, according to the definition of QoS in (\ref{eq:QoS}), charging time is fixed after receiving $f_{i}$, and EV owners can only change the deadline to influence the waiting time.
Therefore, for EV $i$ in the $m$-th charging station, we can also define a virtual queue accounting for the waiting time, which can be written as
\begin{equation} \label{eq:B-virtual-queue}
\!\!\! B_{i, t+1} \!\! = \!\!  \min \! \left\{\! B_{i, t} \!+\!  ({\tt SOC}_{i}^{\tt fin} \!\!-\! {\tt SOC}_{i, t}) \!\left( \pmb{\omega}_{i}^{T} \mathbf{d}_{i} \!-\! f_{i}\right)\!,  B_{i}^{\tt max} \right\}, 
\end{equation}
where $\pmb{\omega}_{i}^{T} \mathbf{d}_{i} - f_{i}$ represents the time difference between the selected deadline and the estimated finishing charging time.
Also, in our design, the waiting time decreases if the SOC approaches the desired SOC, ${\tt SOC}_{i}^{\tt fin}$. 
Then, we restrict the value of the virtual queue with $B_{i}^{\tt max}$.
According to (\ref{eq:B-virtual-queue}), we can define another Lyapunov function as
\begin{equation} \label{eq:second_lyapunov_func}
L_{t}^{\tt dead} = \frac{1}{2} \sum_{i \in {\cal H}_{m, t} }  B_{i, t}^{2} .
\end{equation}
The Lyapunov drift based on $L_{t}^{dead}$ will be
\begin{equation}  \label{eq:second_drift}
\Delta_{t}^{\tt dead} = \mathop{\mathbb{E}} \left[ L_{t+1}^{\tt dead} - L_{t}^{\tt dead} | B_{i, t}, f_{i}, {\tt SOC}_{i, t} \right].
\end{equation}
As in Section \ref{subsec:power_decision_queue}, we use \textit{drift-plus-penalty} algorithm, and therefore we have to find the upper bound of the Lyapunov drift.
The upper bound of the Lyapunov drift, if only one EV staying in the charging station, will be
\begin{equation} \label{eq:dead_queue_upper}
B_{i, t+1}^{2} -  B_{i, t}^{2}  \leq  \left( \pmb{\omega}_{i}^{T} \mathbf{d}_{i} - f_{i} \right)^{2} + 2 B_{i, t} \left( \pmb{\omega}_{i}^{T} \mathbf{d}_{i} - f_{i}\right).
\end{equation}
The detailed derivation is provided in Appendix \ref{subsec:second_queue}.
Based on (\ref{eq:dead_queue_upper}), we can then define 
\begin{equation} \label{eq:J2}
 J_{2}(\mathbf{W}) = \sum_{i \in {\cal H}_{m, t}} \left( \left( \pmb{\omega}_{i}^{T} \mathbf{d}_{i} - f_{i} \right)^{2} + 2 B_{i, t} \left( \pmb{\omega}_{i}^{T} \mathbf{d}_{i} - f_{i}\right) \right).
\end{equation} 
By combining (\ref{eq:second_drift}) with (\ref{eq:J2}), we can get following relation 
\begin{equation}
\Delta_{t}^{\tt dead} \leq \mathop{\mathbb{E}} \left[ J_{2}(\mathbf{W}) | B_{i, t}, f_{i}, {\tt SOC}_{i, t}\right].
\end{equation}
Identical to the steps in Section \ref{subsec:power_decision_queue}, we minimize $J_{2}$ and then add a penalty term followed by \textit{drift-plus-penalty} algorithm.
The problem can now be formulated as
\begin{align}\label{eq:P2}
\mathcal{P}_2:  \min_{\qw} &  ~~ J_{2}(\mathbf{W})  \!+\!\!\!\! \sum_{i \in\! {\cal H}_{m, t} } \sum_{j = 1}^{k} V_{i}^{\tt dead} w_{i, j} ({\tt SOC}_{i}^{\tt fin} \!-\! {\tt SOC}_{i, t}) e^{f_{i} \!-\! d_{i, j}}   \notag\\
\mbox{subject to}    &  \quad (\ref{eq:NE_weight_sum_per_EV}), (\ref{eq:NE_weight_cons}) 
\end{align} 
The objective function of (\ref{eq:P2}) is also a trade-off between changing the deadline and the risk function.
In other words, EV owner $i$ may choose to leave the charging station earlier if $f_{i}$ is earlier than $v_{i}$; however, the risk function goes higher for this choice.
That is, $f_{i}$ may be influenced by the renewable power, such that SOC state may not reach  ${\tt SOC}_{i}^{fin}$ if EV owner $i$ wants to leave earlier.
problem $\mathcal{P}_2$ is a quadratic optimization problem, which can also be solved by using interior-point method \cite{boyd-cvx-book}.

\section{Online Algorithm Design}

In order to search the Nash equilibria for both stochastic games presented in Section \ref{sec:game_model}, we convexify the original formulations to two optimization problems in Section \ref{sec:problem_trans}.
In this section, we design an online algorithm based on these formulations and then provide the theoretical analysis.

\subsection{Finishing Charging Time Estimation} \label{subsec:finish_time_est}

In the deadline selection game, the charging station will broadcast $f_{i}$ for every EV owner.
Hence, before introducing the algorithm in detail, we discuss how the charging station can estimate $f_{i}$.
In this case, we need the aid of the forecasted base load information.
Also, the future electricity price is hard to be forecasted such that it is unknown to the charging stations.
Hence, we can only assume the electricity price is the same in the future.
This assumption is only used for estimating $f_{i}$ for EV owners.
We can then use the following remark to get the idea of estimating the finishing charging time.
\begin{Remark}\label{remark:finish_time_est}
The optimal solution of EV charging problem with the same cost function is to balance the load profile after charging \cite{2015-mpc-EVselect,2018-chung-EV-selection}.
That is, the total charging demand should be distributed to the time slots from $t$ to $t' = \{ \argmin v_{i} | i \in {\cal H}_{t} \}$.
With the based load information and the charging demand, $f_{i}$ for EV owners can be estimated.
\end{Remark}

\subsection{Stochastic EV Charging Management Algorithm}\label{subsec:total_algorithm}

At the beginning of the time slot, each charging station observes the current base load and the current generation of renewable energy information. 
A charging station also determines the set of EVs and the corresponding charging demand in the charging station.
With this information, the charging stations solve Problem ${\cal P}_{1}$ to determine the amount of power to draw from the grid by applying numerical optimization algorithms \cite{boyd-cvx-book}.
The virtual queues $Q_{m, t}$ and $Z_{m, t}$ can be updated with (\ref{eq:Q-virtual-queue}) and (\ref{eq:Z-virtual-queue}) according to the solution of Problem ${\cal P}_{1}$, respectively.
After knowing the available energy for charging EVs, the charging station decides which EV can be charged.
However, in this paper, we do not focus on selecting EVs to be charged, as the literature offers many solutions to this aspect.
Therefore, we employ the idea of the earliest deadline first (EDF) method as in \cite{2015-mpc-EVselect}, where the EDF method assigns the charging rate to EVs according to the deadline.
After receiving power from the charging station, the SOC of EV $i$ can be updated as 
\begin{equation}\label{eq:SOC_update}
{\tt SOC}_{i, t+1} = {\tt SOC}_{i, t} + \epsilon_{c} \frac{P_{i, t}}{E_{i}^{\rm cap}}.
\end{equation}
The charging station then estimates $f_{i}$ for EV owners based on the Remark \ref{remark:finish_time_est} in Section \ref{subsec:finish_time_est}.
With the updated SOC information and $f_{i}$, EV owners can decide to shift their deadline by solving Problem ${\cal P}_{2}$.
The deadline can then be updated using
\begin{equation}\label{eq:deadline_update}
\hat{v_{i}} = \{ d_{i, k} | k = \argmax ~\omega_{i, k}  \}.
\end{equation}
The $B_{i, t}$ can be updated based on the solution of Problem ${\cal P}_{2}$.
Also, the charging station has to collect the charging information of EVs, which arrive at time $t$, at the end of the time slot.
Based on the selection of the EV owners, the $m$-th charging station updates the $V_{m}^{\tt charg}$ as 
\begin{equation} \label{eq:parameter_update}
\!\!\! V_{m}^{\tt charg} = \left\{
\begin{array}{ll}
                \!\!(1 - \alpha) V_{m}^{\tt charg},      & \!\!\mbox{if}  \sum_{i \in {\cal H}_{m, t}} (\hat{v_{i}} - v_{i} )  < 0 ,  \\
                \!\!(1 + \alpha) V_{m}^{\tt charg},      & \!\!\mbox{if}  \sum_{i \in {\cal H}_{m, t}} (\hat{v_{i}} - v_{i} )  > 0 ,  \\
                \!\!V_{m}^{\tt charg},      & \!\!\mbox{if}  \sum_{i \in {\cal H}_{m, t}} (\hat{v_{i}} - v_{i} )  = 0 ,  \\
\end{array}
        \right.
\end{equation}
where $\alpha \approx 0$.
The intuition of (\ref{eq:parameter_update}) is that charging stations decrease the trade-off term to purchase more power from the grid if most of the EV owners wish to leave earlier, and vice versa.
The steps repeat until the time slot reaches ${\cal T}$.
On the basis of the discussion above, we summarize the steps in Algorithm \ref{ago:EV_CSMA}.

\begin{algorithm}
\caption{Stochastic EV Charging Management Algorithm}
\label{ago:EV_CSMA}
 \DontPrintSemicolon
\KwIn{$N$ EVs with their information}
\KwOut{ $x_{m, t}$, ${\tt P}_{i, t}$, $\hat{v_{i}}$}
\For{t  = 1 \KwTo ${\cal T}$ }{
	Determine ${\cal H}_{m, t}$ and ${\tt E}_{i}^{\tt de}$ at time $t$\;
	Get the real-time electricity price information at time $t$\;
	Solve Problem ${\cal P}_{1}$ \\
	Update virtual queue with (\ref{eq:Q-virtual-queue}) and (\ref{eq:Z-virtual-queue}) \\
	Select which EVs to be charged by applying EDF \\
	Use (\ref{eq:SOC_update}) to update SOC information \\
	Estimate $f_{i}$ based on the Remark \ref{remark:finish_time_est} \\
	Solve Problem ${\cal P}_{2}$ and choose $\hat{v_{i}}$ with (\ref{eq:deadline_update}) \\
	Update virtual queue with (\ref{eq:B-virtual-queue})	\\
	Charging stations collect the charging information of upcoming EVs\\
	Update the parameters based on (\ref{eq:parameter_update}) 
   }
\end{algorithm}

\subsection{Performance and Computation Complexity Analysis}

After designing the online algorithm in Section \ref{subsec:total_algorithm}, we analyze the performance of the proposed algorithm in this section.
In Algorithm \ref{ago:EV_CSMA}, ${\cal P}_{1}$ and ${\cal P}_{2}$ are convex optimization problems, and therefore the solution of ${\cal P}_{1}$ and ${\cal P}_{2}$ will converge to the optimal solution by utilizing the numerical algorithms from \cite{boyd-cvx-book}.
However, both (\ref{eq:NE_power_decision_game}) and (\ref{eq:NE_deadline_selection_game}) are stochastic optimization problems, and therefore the optimal solutions of ${\cal P}_{1}$ and ${\cal P}_{2}$ can be the sub-optimal solutions of the stochastic optimization problems.
Hence, we need to show the relation between the solutions of ${\cal P}_{1}$ and ${\cal P}_{2}$, and the solutions of (\ref{eq:NE_power_decision_game}) and (\ref{eq:NE_deadline_selection_game}) , respectively.

We first show the relation between the optimal solution of ${\cal P}_{1}$ and (\ref{eq:NE_power_decision_game}).
Let $\lambda_{i, t}$, $r_{m, t}$, and $k_{t}$ be independent and identically distributed (i.i.d) over time horizon.
Then, there exists a randomized stationary policy such that 
\begin{equation} \label{eq:optimal_cost_form}
\mathop{\mathbb{E}} \left[ k_{t}  \sum_{m=1}^{M} x^{*}_{m, t}  \right] = c_{t}^{*},
\end{equation}
where $c_{t}^{*}$ is the optimal cost that can be achieved over time slot $t$ for (\ref{eq:NE_power_decision_game}), and $x^{*}_{m, t}$ is the corresponding optimal control decisions.
\begin{Theorem}\label{theorem:optimal_gap}
The total electricity cost in a day obtained by solving Problem ${\cal P}_{1}$ is within ${\cal O}(\frac{1}{ V_{m}^{\tt charg} })$ to the NE of the cost minimization game, which implies
\begin{equation} 
 \sum_{t = 1}^{ \cal T }  c_{t}^{*}  \leq \sum_{t = 1}^{ \cal T }  \mathop{\mathbb{E}} \left[k_{t} g(\mathbf{x}) \right]  \leq   \frac{D {\cal T} }{ V_{m}^{\tt charg} } + \sum_{t = 1}^{ \cal T }  c_{t}^{*}.
\end{equation}
The proof is provided in Appendix \ref{subsec:proof_theorem_1}.
\end{Theorem}
The relation between the solution of Problem ${\cal P}_{2}$ and NE of deadline selection game can be obtained in a similar manner in Theorem \ref{theorem:optimal_gap}.
With Theorem \ref{theorem:optimal_gap}, we can conclude that the solutions of ${\cal P}_{1}$ and ${\cal P}_{2}$ are the $\epsilon$-optimal solutions of (\ref{eq:NE_power_decision_game}) and (\ref{eq:NE_deadline_selection_game}), respectively.
Therefore, Theorem \ref{theorem:optimal_gap} implies that the proposed algorithm can yield a solution close to the Nash equilibria of the cost minimization and deadline selection games by increasing the trade-off parameters.
However, increasing $V_{m}^{\tt charg}$ represents that the charging stations will focus more on minimizing the electricity cost.
Then, increasing $V_{i}^{\tt dead}$ increases the waiting time for EV owners.

Next, we analyze the computation complexity of Algorithm \ref{ago:EV_CSMA}.
Operations other than those in lines $3$ and $8$ are all linear.
We therefore focus on analyzing the computation complexity of solving the optimization formulation in lines $3$ and $8$.
In lines $3$ and $8$, both Problem ${\cal P}_{1}$ and ${\cal P}_{2}$ are quadratic convex optimization problems, and they can be easily solved.
For the traditional interior-point method, the complexity is ${\cal O}(n \log(n))$, where $n$ is the number of the variables \cite{2004-SIAM-complexity}.
Therefore, the computation complexity of solving Problem ${\cal P}_{1}$ and ${\cal P}_{2}$ are ${\cal O}(M \log(M))$ and ${\cal O}(N \log(N))$, respectively.
The computation complexity of Algorithm \ref{ago:EV_CSMA} will then be ${\cal O}(n \log(n))$, where $n\!=\!\max \{ M, N \}$.
Thus, we can claim that the overall complexity of Algorithm \ref{ago:EV_CSMA} is reasonably low.

\section{Numerical Results}\label{sec:simulation}

We use real-world data from CAISO \cite{Load_data} to evaluate the performance of the proposed algorithm.
The data from 05/01/2018 is used.
The base load information is normalized by the maximum value.
Then, we can get the solar and wind energy capacity in California, such that we can calculate the generation profile of unit renewable energy capacity.
We also employ real-time pricing data from \cite{Load_data}.
Unless otherwise specified, the simulation settings are as follows.

The total number of EVs is $N$, and they are randomly assigned to one of the $M$ charging stations, each of them having $30$ kWh solar and $10$ kWh wind generation capacity.
The set of variables for different scenarios are provided in Table \ref{tb:charging_scenario_setting}.
The base load information is generated by multiplying the maximum load in a day and the normalized load curve.

The EVs used in the simulation are all Nissan Leaf 2018, each with a battery capacity of $40$ kWh and a maximum charging rate of $6.6$ kW.
The initial SOC values are randomly and uniformly generated from the interval $[0,1]$, and the EVs are assumed to be fully charged at the departure time.
The charging efficiency is set to $90\%$. 
The arrival time and deadlines are generated randomly around $8:00$ and $17:00$, respectively.

The time horizon is divided into $96$ time slots with a length of $15$ minutes to represent a $24$-hour period.
Let $V_{m}^{\tt charg}$ be $350$ for all charging stations, $V_{i}^{dead}$ be $200$ for all EV owners, $\eta_{m}$ be $1$, and $\alpha$ be $0.001$ as initial values. 
In the deadline selection game, EV owners can choose between $5$ actions, $\{ v_{i}-2, v_{i}-1, v_{i}, v_{i}+1, v_{i}+2 \}$, which means $k=5$.
Then, $B_{i}^{\tt max}$ is set to $20$.

\begin{table} \small
\begin{center}
\caption{The parameter setting for different scenarios}\label{tb:charging_scenario_setting}
\begin{tabular}{|C{3.4cm}|C{0.6cm}|C{0.6cm}|C{0.6cm}|C{0.6cm}|c|c|c|}
\hline
$N$      					& $25$  & $50$   & $100$	& $200$  \\
\cline{1-5}
$M$ 						& $2$   & $4$  	 & $8$ 		& $16$   \\
\cline{1-5}
Maximum Load (Kw)   		& $100$ & $200$  & $400$  	& $800$ \\
\cline{1-5}	
${\tt P}^{\tt peak} $ (Kw) 	& $120$ & $240$  & $440$    & $850$  \\
\hline
\end{tabular}
\end{center}
\end{table}

During the simulation, we compare our algorithm with the offline charging cost minimization algorithm (OCCMA) in \cite{2012-He-EV-selection} and the online charging scheduling algorithm (OCSA) in \cite{2018-chung-EV-selection} with the user convenience defined in \cite{2012-wen-EV-selection} as the benchmarks.
The OCCMA obtains the real electricity price and the real generation of renewable energy, and therefore the OCCMA leads to the minimum electricity cost.
Then, we analyze two scenarios to OCSA denoted as OCSA-F and OCSA-N.
The forecasted renewable generation and forecasted electricity price are available to the OCSA denoted as OCSA-F.
By contrast, OCSA-N only obtains the information of the forecasted renewable generation.
Forecasted electricity price is generated by adding Gaussian noise (with zero mean and dime variance) to the real-time pricing data.
Comparison between the OCSA-F and the OCSA-N provides some useful insight on the influence of the forecasted electricity price to the OCSA. 
Also, OCSA-N is similar to the methods in \cite{2012-wen-EV-selection,2017-mal-EV-selection}.

\subsection{Electricity Cost}\label{subsec:cost_compare}

\begin{figure}
\begin{center}
\resizebox{3.2in}{!}{%
\includegraphics*{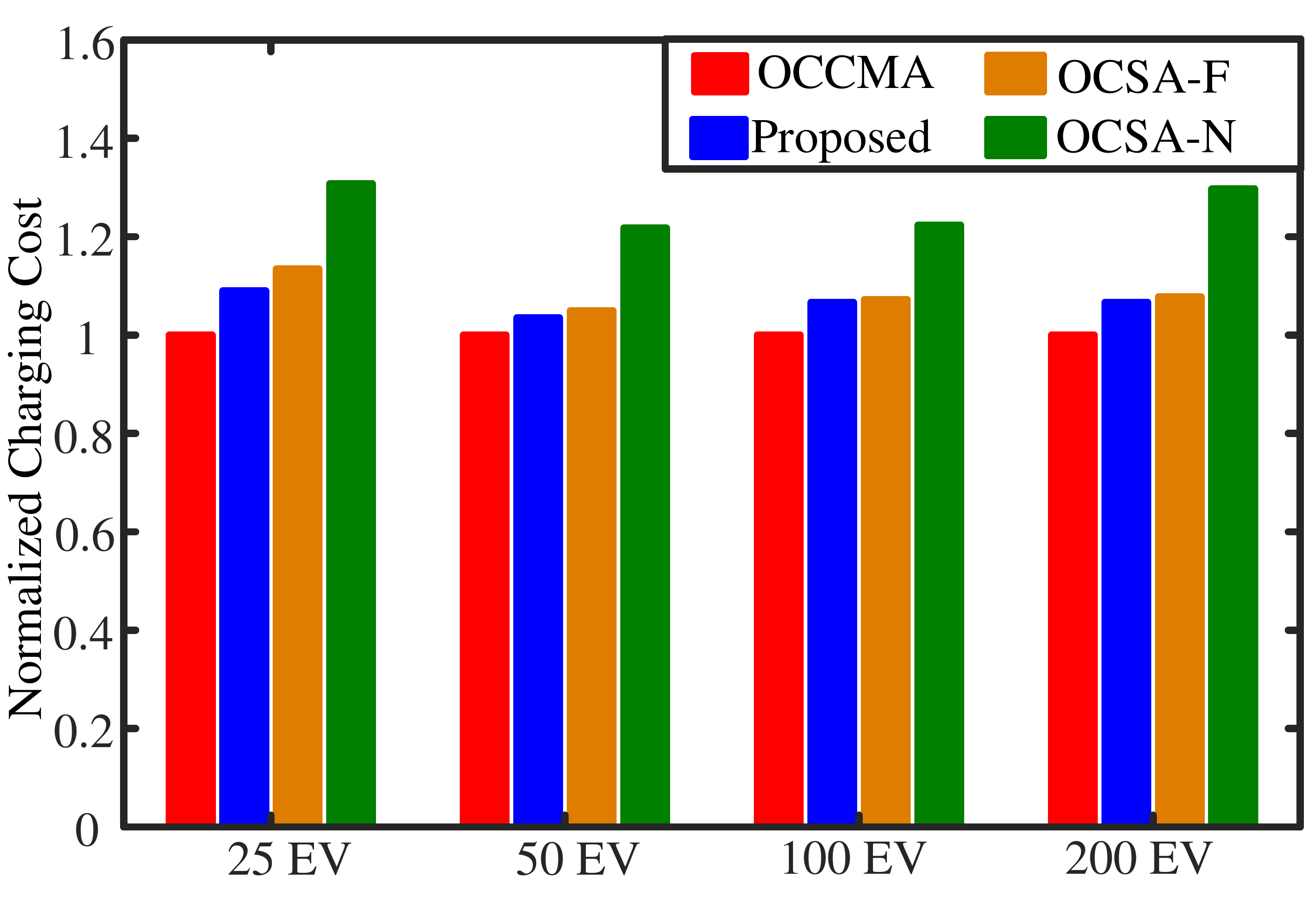} }%
\caption{Charging cost comparison.} 
\label{fig:cost_compare}
\end{center}
\end{figure}

In this part, we consider the electricity cost of different scenarios listed in Table \ref{tb:charging_scenario_setting}.
Then, the results are presented in Fig.\,\ref{fig:cost_compare}.
Note that we normalized the electricity cost with the result of the OCCMA.

According to the simulation results, the proposed algorithm yields results very close to the optimal solution.
The proposed method leads to $7.07\%$ higher cost compared to the optimal solution.
The OCSA-F yields a comparable performance as the proposed method.
However, the difference is $26.12\%$ for OCSA-N. 
By applying the proposed algorithm, we can save about $20.00\%$ of the electricity cost compared to OCSA-N.

The OCSA-F leads to better performance than the OCSA-N because of the knowledge about the future electricity price.
However, the performance of OCSA-F depends heavily on the accuracy of the forecasted electricity price.
If the difference between the forecasted electricity price and the real electricity price increases, the performance of the OCSA-F may degrade to the performance of OCSA-N.
OCSA-N performs poorly because the future electricity price is unknown to the charging stations.
Also, the estimated wind power generation on 05/01/2018 is much lower than the real power generation.
Therefore, the charging stations purchase more power from the grid at the current time slot without considering the electricity price.
By contrast, with the proposed formulation, we have a trade-off term between the remaining charging demand and the electricity cost.
If the price is too high, the charging stations therefore use the renewable energy to charge without purchasing the external energy from the power grid.

\subsection{QoS}

Other than the electricity cost, we study the outcome of the deadline selection game.
First, we present the results of the deadline selection game.
We assume there is only one EV in the charging station for ease of explanation.
Let $f_{i}$ be $16$, $V_{2}$ be $250$, and $B_{i, t}$ be $10$.
Three different SOC, which are $0.52$ (Case 1), $0.64$ (Case 2), and $0.73$ (Case 3), are compared.
We then compare the influence of the value of virtual queue $B_{i, t}$.
Let $f_{i}$ be $16$, $V_{2}$ be $250$, and ${\tt SOC}_{i, t}$ be $0.64$.
Three different virtual queue values, namely $10$ (Case 2), $6.5$ (Case 4), and $3.5$ (Case 5), are discussed.
The EV has five choices, namely $17:30, 17:45, 18:00, 18:15$, and $18:30$.
The results are provided in Table \ref{tb:dead_game_result}.

From the results in Table \ref{tb:dead_game_result}, it is clear that EV owners postpone the deadline if they have low SOC right now. 
By contrast, the EV owner wishes to leave the charging station earlier if they have higher SOC at the current time slot.
Thus, EV owners have higher confidence to leave earlier if they are at higher SOC.
Next, we discuss the results of changing the value of $B_{i, t}$.
Under the same SOC, the EV owner with a higher value of $B_{i, t}$ would like to prepone the deadline.
That is, if the EV owners stay a long time in the charging station without receiving power, the EV owner tends to leave earlier.

\begin{table} \small
\begin{center}
\caption{The probability of selecting the deadline}
\label{tb:dead_game_result}
\begin{tabular}{|c|c|c|c|c|c|}
\hline
           		& $17:30$ & $17:45$  & $18:00$   & $18:15$		& $18:30$ \\
\hline \hline
Case 1      	& $0$  	  & $0$   	 & $0.2610$  & $0.7390$		& $0$ \\
\hline
Case 2  		& $0$  	  & $0.3246$ & $0.6754$  & $0$			& $0$ \\
\hline
Case 3	   		&$0.3483$ & $0.6517$ & $0$		 & $0$			& $0$ \\
\hline \hline
Case 2  		& $0$  	  & $0.3246$ & $0.6754$  & $0$			& $0$ \\
\hline
Case 4			& $0$  	  & $0$ 	 & $0.4466$  & $0.5534$		& $0$ \\
\hline
Case 5  		& $0$  	  & $0$ 	 & $0$  	 & $0.2154$		& $0.7846$ \\
\hline
\end{tabular}
\vspace{-0.3cm}
\end{center}
\end{table}

\begin{figure}
\begin{center}
\resizebox{3.2in}{!}{%
\includegraphics*{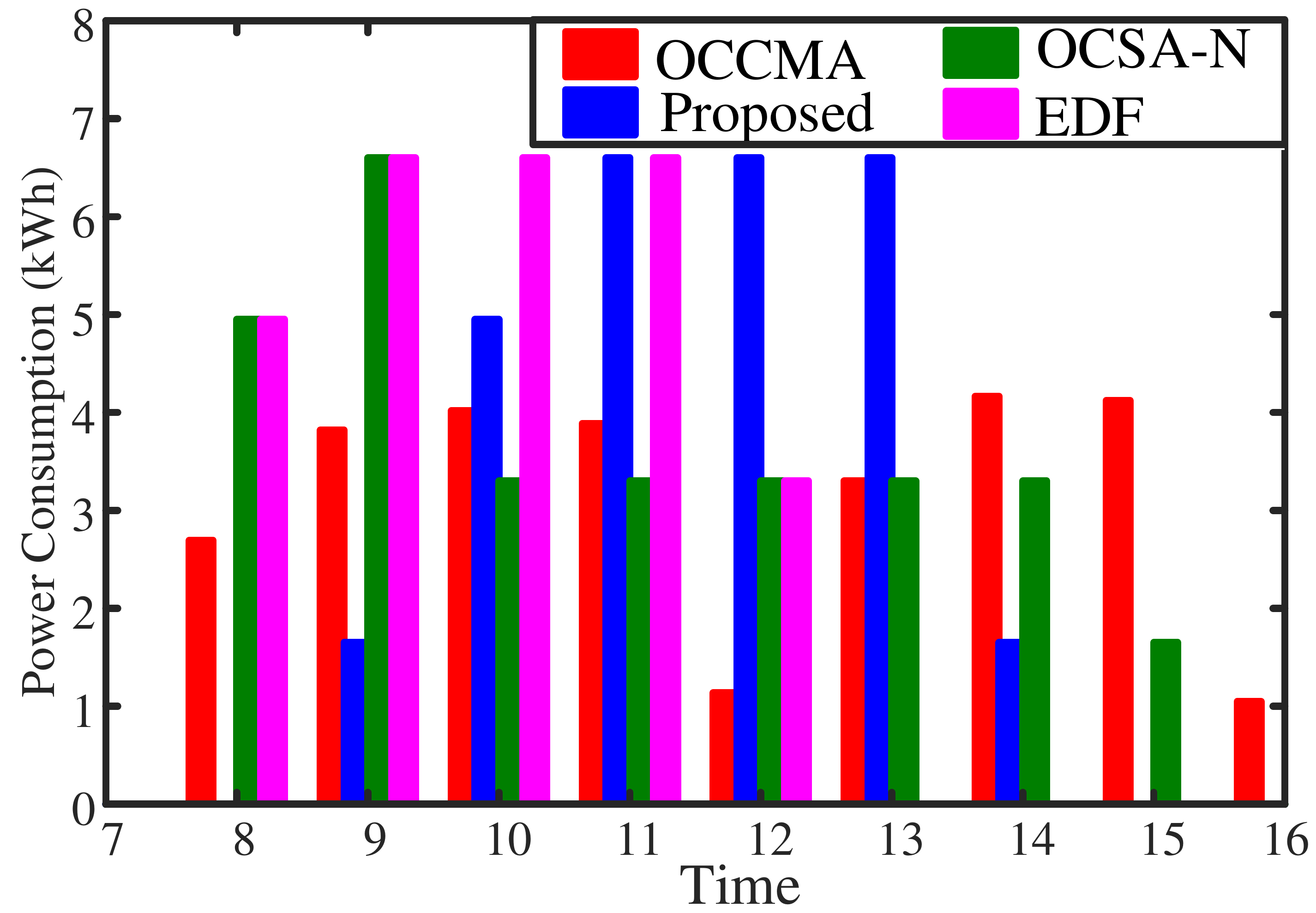} }%
\caption{Charging profile of the $14$-th EV.} 
\label{fig:profile_compare_14}
\end{center}
\end{figure}

After showing the outcome of the deadline selection game, we discuss its influence on the charging process.
Therefore, we need to detailly analyze the charging profile of the EV.
The charging profile of the $14$-th EV is provided in Fig.\ref{fig:profile_compare_14}.
Here, we also compare the proposed method without deadline selection game, which is EDF.
OSCA-F has a similar charging mechanism as OCSA-N, and thus the charging profile of OSCA-F is omitted. 
The charging time of the $14$-th EV is $8.75$, $7.75$, $6.75$, and $5.75$ hrs by utilizing OCCMA, OCSA-N, proposed method, and EDF, respectively.	
We can observe that the charging time of the proposed method is higher than EDF but lower than OCCMA and OCSA-N.
For both OCCMA and OCSA-N, the EV keeps receiving power until the deadline is reached; therefore, the charging time is long. 
In the EDF method, the $14$-th EV will be immediately charged after it enters the charging station because it has an earlier deadline than most of the EVs.
However, for the proposed method, it will postpone the charging task around $1$ hour because some EVs already have accumulative waiting time.
Hence, the $14$-th EV needs to wait to receive the charging power.

\begin{figure}
\begin{center}
\resizebox{3.2in}{!}{%
\includegraphics*{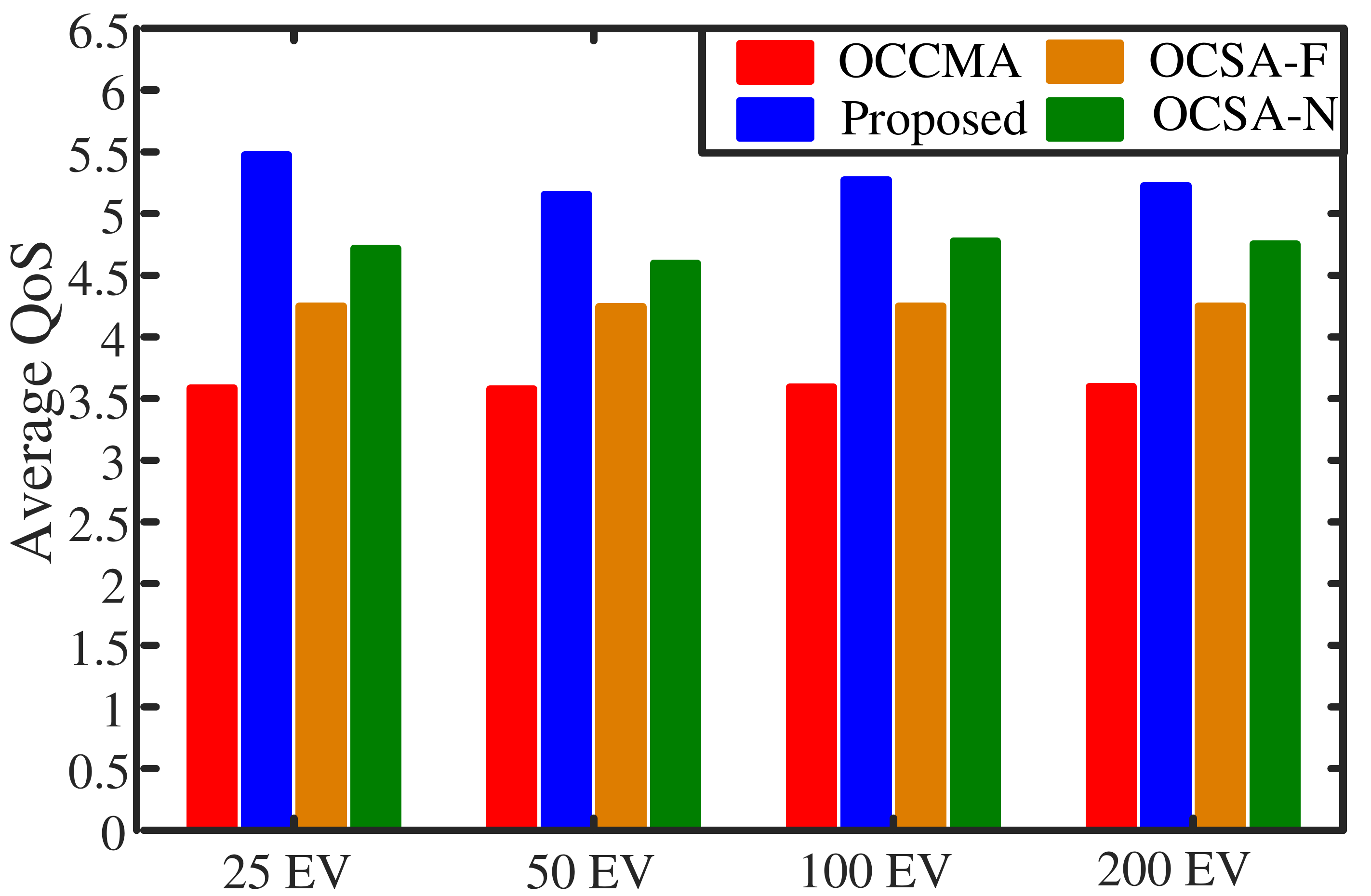} }%
\caption{Average QoS evaluation.} 
\label{fig:QoS_compare}
\end{center}
\end{figure}

We then evaluate the average QoS for the EV owners in the charging stations in Fig.\,\ref{fig:QoS_compare}.
According to the results, the proposed algorithm can obtain $47.19\%$ and $12.13\%$ higher QoS than the benchmarks, respectively.
For the OCCMA, it distributes the charging demands between $a_{i}$ and $v_{i}$, such that the charging tasks are finished just before the deadline.
In this case, the charging time is long, but the waiting is short.
By contrast, user satisfaction is used to schedule the charging tasks without considering the cost in OCSA-N.
By using this algorithm, some EVs constantly receive the power, but some need to wait a long time to get the power from the charging station.
The charging time of the OCSA-N thus varies considerably.
Then, the average QoS of the OCSA-F lies between the OCCMA and the OCSA-N because of the accuracy of the forecasted information as mentioned in Sec. \ref{subsec:cost_compare}.
Therefore, some EVs obtain high QoS but some EVs have low QoS; however, the average QoS is still higher than the OCCMA.
If the forecasted information is accurate, the average QoS of the OCSA-F will be the same as OCCMA.
On the other hand, the average QoS of the OCSA-F can approach to the average QoS the OCSA-N if the forecasted information is inaccurate.
For the proposed method, we allow EVs to shift their deadline, such that we can avoid the issue of a large variance in charging time in the OSCA-N.
That is why our proposed method can obtain higher average QoS.

\section{Conclusion}\label{sec:conclusion}

In this work, we presented an intelligent EV charging management scheme for the charging stations in the transportation system.
By using the proposed method, we can eliminate the negative impact of charging a large and rapidly growing number of EVs in electric mobility to the power grid.
We modeled the interactions between the power grid and charging stations as a stochastic game.
We incorporated the dynamic behavior of EV owners in terms of their preference or choice on charging parameters as another stochastic game.
We then proposed an online algorithm to approach the Nash equilibria of both games.
We utilized real data from California Independent System Operator to evaluate the performance of the proposed algorithm.
The numerical results illustrated that the proposed algorithm can achieve electricity cost very close to the minimum electricity cost while also enhancing QoS significantly.

{\renewcommand{\baselinestretch}{1}
\begin{footnotesize}
\bibliographystyle{IEEEtran}
\bibliography{References_TITS}

\begin{thebibliography}{10}
\providecommand{\url}[1]{#1}
\csname url@samestyle\endcsname
\providecommand{\newblock}{\relax}
\providecommand{\bibinfo}[2]{#2}
\providecommand{\BIBentrySTDinterwordspacing}{\spaceskip=0pt\relax}
\providecommand{\BIBentryALTinterwordstretchfactor}{4}
\providecommand{\BIBentryALTinterwordspacing}{\spaceskip=\fontdimen2\font plus
\BIBentryALTinterwordstretchfactor\fontdimen3\font minus
  \fontdimen4\font\relax}
\providecommand{\BIBforeignlanguage}[2]{{%
\expandafter\ifx\csname l@#1\endcsname\relax
\typeout{** WARNING: IEEEtran.bst: No hyphenation pattern has been}%
\typeout{** loaded for the language `#1'. Using the pattern for}%
\typeout{** the default language instead.}%
\else
\language=\csname l@#1\endcsname
\fi
#2}}
\providecommand{\BIBdecl}{\relax}
\BIBdecl

\bibitem{2017-EV-outlook}
\BIBentryALTinterwordspacing
{IEA Publications}, ``{Global EV Outlook 2017: Two Million and Counting},''
  International Energy Agency, Tech. Rep., 2017. [Online]. Available:
  \url{https://webstore.iea.org/global-ev-outlook-2017}
\BIBentrySTDinterwordspacing

\bibitem{2018-EV-outlook}
\BIBentryALTinterwordspacing
------, ``{Global EV Outlook 2018: Towards Cross-modal Electrification},''
  International Energy Agency, Tech. Rep., 2018. [Online]. Available:
  \url{https://webstore.iea.org/global-ev-outlook-2018}
\BIBentrySTDinterwordspacing

\bibitem{2015-ev-impact}
{A. Dubey and S. Santoso}, ``{Electric Vehicle Charging on Residential
  Distribution Systems: Impacts and Mitigations},'' \emph{IEEE Access}, vol.~3,
  pp. 1871--1893, Sep. 2015.

\bibitem{2012-He-EV-selection}
{Y. He, B. Venkatesh, and L. Guan}, ``{Optimal Scheduling for Charging and
  Discharging of Electric Vehicles},'' \emph{IEEE Trans. Smart Grid}, vol.~3,
  no.~3, pp. 1095--1105, Sep. 2012.

\bibitem{2018-shi-EV-selec}
{Y. Shi, H. D. Tuan, A. V. Savkin, T. Q. Duong, and H. V. Poor}, ``{Model
  Predictive Control for Smart Grids with Multiple Electric-Vehicle Charging
  Stations},'' \emph{IEEE Trans. Smart Grid}, vol.~10, no.~2, pp. 2127--2136,
  Mar. 2019.

\bibitem{2018-david-EV-manage}
{C. P. Mediwaththe and D. B. Smith}, ``{Game-Theoretic Electric Vehicle
  Charging Management Resilient to Non-Ideal User Behavior},'' \emph{IEEE
  Trans. Intell. Transp. Syst.}, vol.~19, no.~11, pp. 3486--3495, Nov. 2018.

\bibitem{2019-Seyedyazdi-EV-manage}
{M. Seyedyazdi, M. Mohammadi, and E. Farjah}, ``{A Combined Driver-Station
  Interactive Algorithm for a Maximum Mutual Interest in Charging Market},''
  \emph{IEEE Trans. Intell. Transp. Syst.}, 2019, to be published.

\bibitem{2019-zhang-EV-manage}
{Y. Zhang, P. You, and L. Cai}, ``{Optimal Charging Scheduling by Pricing for
  EV Charging Station With Dual Charging Modes},'' \emph{IEEE Trans. Intell.
  Transp. Syst.}, vol.~20, no.~9, pp. 3386--3396, Sep. 2019.

\bibitem{2012-wen-EV-selection}
{C.-K. Wen, J.-C. Chen, J.-H. Teng, and P.-A. Ting}, ``{Decentralized Plug-in
  Electric Vehicle Charging Selection Algorithm in Power Systems},'' \emph{IEEE
  Trans. Smart Grid}, vol.~3, no.~4, pp. 1779--1789, Dec. 2012.

\bibitem{2017-mal-EV-selection}
{A. Malhotra, G. Binetti, A. Davoudi, and I. D. Schizas}, ``{Distributed Power
  Profile Tracking for Heterogeneous Charging of Electric Vehicles},''
  \emph{IEEE Trans. Smart Grid}, vol.~8, no.~5, pp. 2090--2099, Sep. 2017.

\bibitem{2018-chung-EV-selection}
{H.-M Chung, W.-T. Li, C. Yuen, C.-K Wen, and N. Crespi}, ``{Electric Vehicle
  Charge Scheduling Mechanism to Maximize Cost Efficiency and User
  Convenience},'' \emph{IEEE Trans. Smart Grid}, vol.~10, no.~3, pp.
  3020--3030, May 2019.

\bibitem{2015-zhao-stochastic-infocom}
{S. Zhao, X. Lin, and M. Chen}, ``{Peak-minimizing Online EV Charging:
  Price-of-uncertainty and Algorithm Robustification},'' in \emph{Proc. IEEE
  Intl. Conf. on Comput. Commun. (INFOCOM)}, Kowloon, Hong Kong, April 2015,
  pp. 2335--2343.

\bibitem{2016-wayes-EV-manage}
{W. Tushar, C. Yuen, S. Huang, D. B. Smith, and H. V. Poor}, ``{Cost
  Minimization of Charging Stations With Photovoltaics: An Approach With EV
  Classification},'' \emph{IEEE Trans. Intell. Transp. Syst.}, vol.~17, no.~1,
  pp. 156--169, Jan. 2016.

\bibitem{2018-luo-stochastic}
{C. Luo, Y. F. Huang, and V. Gupta}, ``{Stochastic Dynamic Pricing for EV
  Charging Stations With Renewable Integration and Energy Storage},''
  \emph{IEEE Trans. Smart Grid}, vol.~9, no.~2, pp. 1494--1505, Mar. 2018.

\bibitem{2017-huang-stochastic}
{Q. Huang, Q. S. Jia, and X. Guan}, ``{A Multi-Timescale and Bilevel
  Coordination Approach for Matching Uncertain Wind Supply With EV Charging
  Demand},'' \emph{IEEE Trans. Autom. Sci. Eng.}, vol.~14, no.~2, pp. 694--704,
  Apr. 2017.

\bibitem{2018-huang-stochastic}
------, ``{Robust Scheduling of EV Charging Load With Uncertain Wind Power
  Integration},'' \emph{IEEE Trans. Smart Grid}, vol.~9, no.~2, pp. 1043--1054,
  Mar. 2018.

\bibitem{2018-yang-stochastic}
{Y. Yang, Q.-S. Jia, G. Deconinck, X. Guan, Z. Qiu, and Z. Hu}, ``{Distributed
  Coordination of EV Charging With Renewable Energy in a Microgrid of
  Buildings},'' \emph{IEEE Trans. Smart Grid}, vol.~9, no.~6, pp. 6253--6264,
  Nov. 2018.

\bibitem{2018-hiroshi-stochastic}
{H. Kikusato, K. Mori, S. Yoshizawa, Y. Fujimoto, H. Asano, Y. Hayashi, A.
  Kawashima, S. Inagaki, and T. Suzuki}, ``{Electric Vehicle Charge-Discharge
  Management for Utilization of Photovoltaic by Coordination between Home and
  Grid Energy Management Systems},'' \emph{IEEE Trans. Smart Grid}, vol.~10,
  no.~3, pp. 3186--3197, May 2019.

\bibitem{2018-zhou-EV-lyapunov}
{Y. Zhou, D. K. Y. Yau, P. You, and P. Cheng}, ``{Optimal-Cost Scheduling of
  Electrical Vehicle Charging Under Uncertainty},'' \emph{IEEE Trans. Smart
  Grid}, vol.~9, no.~5, pp. 4547--4554, Feb. 2018.

\bibitem{2018-gigoni-forecast-error}
{L. Gigoni, A. Betti, E. Crisostomi, A. Franco, M. Tucci, F. Bizzarri, and D.
  Mucci}, ``{Day-Ahead Hourly Forecasting of Power Generation From Photovoltaic
  Plants},'' \emph{IEEE Trans. Sustain. Energy}, vol.~9, no.~2, pp. 831--842,
  Apr. 2018.

\bibitem{2016-wang-DR-user}
{Y. Wang, W. Saad, N. B. Mandayam, and H. V. Poor}, ``{Load Shifting in the
  Smart Grid: To Participate or Not?}'' \emph{IEEE Trans. Smart Grid}, vol.~7,
  no.~6, pp. 2604--2614, Nov. 2016.

\bibitem{2014-drive-range}
{S. A. Birrell, A. McGordon, and P. A. Jennings}, ``{Defining the Accuracy of
  Real-World Range Estimations of an Electric Vehicle},'' in \emph{Proc. 17th
  Intl. IEEE Conf. Intell. Transp. Syst. (ITSC)}, Qingdao, China, Oct. 2014.

\bibitem{Load_data}
\BIBentryALTinterwordspacing
{CAISO}, ``{California Independent System Operator Open Access Same-time
  Information System}.'' [Online]. Available:
  \url{http://oasis.caiso.com/mrioasis/logon.do}
\BIBentrySTDinterwordspacing

\bibitem{1998-rate-allocation}
{F. P. Kelly, A. K. Maulloo, and D. K. H. Tan}, ``{Rate Control for
  Communication Networks: Shadow Prices, Proportional Fairness and
  Stability},'' \emph{J. Oper. Res. Soc.}, vol.~49, no.~3, pp. 237--252, Apr.
  1998.

\bibitem{2010-michael-lyapunov}
{M. J. Neely}, ``{Universal Scheduling for Networks with Arbitrary Traffic,
  Channels, and Mobility},'' in \emph{Proc. IEEE Conf. on Decision and Control
  (CDC)}, Atlanta, GA, USA, Dec. 2010.

\bibitem{2010-michael-book}
------, ``{Stochastic Network Optimization with Application to Communication
  and Queueing Systems},'' \emph{{Synthesis Lectures on Commun. Netw.}},
  vol.~3, no.~1, pp. 1--211, 2010.

\bibitem{boyd-cvx-book}
S.~Boyd and L.~Vandenberghe, \emph{{Convex Optimization}}.\hskip 1em plus 0.5em
  minus 0.4em\relax New York, NY, USA: Cambridge University Press, 2004.

\bibitem{2015-mpc-EVselect}
{W. Tang and Y.-J. Zhang}, ``{A Model Predictive Control Approach for
  Low-Complexity Electric Vehicle Charging Scheduling: Optimality and
  Scalability},'' \emph{IEEE Trans. Power Syst.}, vol.~32, no.~2, pp.
  1050--1063, Mar. 2017.

\bibitem{2004-SIAM-complexity}
{Y. Q. Bai, M. El Ghami, and C. Roos}, ``{A Comparative Study of Kernel
  Functions for Primal-Dual Interior-Point Algorithms in Linear
  Optimization},'' \emph{SIAM J. Optim.}, vol.~15, no.~1, pp. 101--128, 2004.

\bibitem{2012-michael-lyapunov}
{M. J. Neely}, ``{Stability and Probability 1 Convergence for Queueing Networks
  via Lyapunov Optimization},'' \emph{J. Appl. Math.}, 2012.

\bibitem{2008-michael-lyapunov}
{M. J. Neely and R. Urgaonkar}, ``{Opportunism, Backpressure, and Stochastic
  Optimization with the Wireless Broadcast Advantage},'' in \emph{Proc. 42nd
  Asilomar conf. Signal, Syst., and Comput.}, Pacific Grove, California, USA,
  Oct. 2008.

\end{thebibliography}
\end{footnotesize}}

\appendix

\subsection{Lyapunov Drift upper Bound in Cost Minimization Game} \label{subsec:first_queue}

In (\ref{eq:first_lyapunov_func}), it contains two components, $Q_{m, t}$ and $Z_{m, t}$.
To get the upper bound of (\ref{eq:charging_drift}), we derivate the upper bound of the squared difference of each component separately.
The detailed derivation of the upper bound of $Q_{m, t+1}^{2} - Q_{m, t}^{2}$ is provided as follow:
\begin{subequations} \label{eq:Q_drift_upper}
\begin{align}
& Q_{m, t+1}^{2} = \left(\max\{Q_{m, t} - \epsilon_{c} Y_{m, t}, 0\} + \lambda_{m, t} \right)^{2} \label{eq:square_Q}\\ 
& \leq \left( Q_{m, t} - \epsilon_{c} Y_{m, t} \right)^{2} + 2 \lambda_{m, t} \left( Q_{m, t} - \epsilon_{c} Y_{m, t} \right) +  \lambda_{m, t}^{2}  \label{eq:first_ine_Q}\\
& \leq  Q_{m, t}^{2} + \epsilon_{c}^{2} Y_{m, t}^{2} +  \lambda_{m, t}^{2}  + 2 Q_{m, t} \left( \lambda_{m, t} - \epsilon_{c} Y_{m, t} \right)   \label{eq:third_ine_Q}\\
& \leq  Q_{m, t}^{2} + \epsilon_{c}^{2} Y_{m, t}^{2}  +  \lambda_{m, t}^{\tt max^{2}}  + 2 Q_{m, t} \left( \lambda_{m, t}^{\tt max} - \epsilon_{c} Y_{m, t} \right) \label{eq:final_ine_Q}
\end{align}
\end{subequations}
In (\ref{eq:square_Q}), we square the both side of (\ref{eq:Q-virtual-queue}).
Then, (\ref{eq:first_ine_Q}) is obtained by applying binomial theorem and $\max\{a, 0 \}^{2} \leq a^{2}$.
By neglecting the $2\lambda_{m, t} Y_{m, t}$, we can get (\ref{eq:third_ine_Q}).
In the final, we assign $\lambda_{m, t}$ with $\lambda_{m, t}^{\tt max}$ to get (\ref{eq:final_ine_Q}).

Then, a similar calculation can be carried out for the derivation process of getting the upper bound of $Z_{m, t+1}^{2} - Z_{m, t}^{2}$.
The detailed derivation of the upper bound is provided as follow:
\begin{subequations}
\begin{align}
& \!\!\!\! Z_{m, t+1}^{2} = \left( \max \{Z_{m, t} +  \eta_{m} Q_{m, t} - \epsilon_{c} Y_{m, t} , 0 \} \right)^{2} \\
& \!\!\!\! \leq \left( Z_{m, t} +  \eta_{m} Q_{m, t} - \epsilon_{c} Y_{m, t} \right)^{2} \\
& \!\!\!\! =  Z_{m, t}^{2} \!+\! 2 Z_{m, t} \left(\eta_{m} Q_{m, t} \!-\! \epsilon_{c} Y_{m, t} \right)  +  ( \eta_{m} Q_{m, t} \!-\! \epsilon_{c} Y_{m, t} )^{2} 
\end{align}
\end{subequations}

\subsection{Lyapunov Drift upper Bound in Deadline Selection Game} \label{subsec:second_queue}

In this appendix, we aim to provide the derivation of the upper bound of (\ref{eq:second_drift}).
Eq.\,(\ref{eq:second_lyapunov_func}) contains the summation of all EVs in the charging station.
However, it is easier to show derivation with only one EV in the charging station. Therefore, we can ignore the summation first.
The detailed derivation of the upper bound is provided as follow:
\begin{subequations}
\begin{align}
& B_{i, t+1}^{2} \notag \\
& = \! \left( \!\min \!\left\{ \!B_{i, t}  \!+\!  ( {\tt SOC}_{i}^{\tt fin} \!-\! {\tt SOC}_{i, t} \!) \!\left( \pmb{\omega}_{i}^{T} \mathbf{d}_{i} \!-\! f_{j} \right)\!, B_{max} \right\}  \right)^{2}   \\
&\leq \left( B_{i, t} +  ({\tt SOC}_{i}^{\tt fin}- {\tt SOC}_{i, t} ) \left( \pmb{\omega}_{i}^{T} \mathbf{d}_{i} - f_{j} \right) \right)^{2}  \label{eq:first_ine_B}\\
& =  B_{i, t}^{2} + ({\tt SOC}_{i}^{\tt fin} - {\tt SOC}_{i, t})^{2} \left( \pmb{\omega}_{i}^{T} \mathbf{d}_{i} - f_{i}\right)^{2}  \notag \\
& ~~~+ 2 B_{i, t} ( {\tt SOC}_{i}^{\tt fin} - {\tt SOC}_{i, t}) \left(  \pmb{\omega}_{i}^{T} \mathbf{d}_{i} - f_{i}\right)  \label{eq:third_eq_B}\\
& \leq B_{i, t}^{2} + \left( \pmb{\omega}_{i}^{T} \mathbf{d}_{i} - f_{i} \right)^{2} + 2 B_{i, t} \left( \pmb{\omega}_{i}^{T} \mathbf{d}_{i} - f_{i}\right)  \label{eq:final_ine_B}
\end{align}
\end{subequations}
We still square both side of (\ref{eq:B-virtual-queue}) first. 
Then, the first inequality, (\ref{eq:first_ine_B}), comes from the relation of $\min\{a, b \}^{2} \leq a^{2}$.
We use binomial theorem to get (\ref{eq:third_eq_B}).
In the final, (\ref{eq:final_ine_B}) is obtained by assuming ${\tt SOC}_{i}^{\tt fin} - {\tt SOC}_{i, t} \leq 1 - {\tt SOC}_{i, t}$.

\subsection{Proof of Theorem 1}\label{subsec:proof_theorem_1}

Problem ${\cal P}_{1}$ is designed to minimize the drift bound, which holds for all policies, including the optimal policy given in (\ref{eq:optimal_cost_form}).
For notational simplicity, we assume $V_{m}^{\tt charg}$ is the same for all charging stations.
Let us define $g(\mathbf{x}) = \sum_{m=1}^{M} x_{m, t}$. 
We then introduce the optimal solution to the right-hand side of \textit{drift-plus-penalty} term \cite{2012-michael-lyapunov,2008-michael-lyapunov} as
\begin{align} \label{eq:drift_with_optimal}
& L_{t+1}^{\tt charg} - L_{t}^{\tt charg} + V_{m}^{\tt charg} ~ k_{t} g(\mathbf{x})  \leq    \notag \\
& \frac{1}{2} \sum_{m=1}^{M}  \lambda_{m, t}^{\tt max^{2}} +   J_{1} (\mathbf{x}^{*} )_{| Q_{m, t}, Y_{m, t}, \lambda_{m, t} }  + c_{t}^{*} V_{m}^{\tt charg}.
\end{align}
Taking expectation on both side of (\ref{eq:drift_with_optimal}) and using iterated expectation, we can get following relation
\begin{equation} 
\!\!\! \mathop{\mathbb{E}} \left[ L_{t+1}^{\tt charg} \right]  -  \mathop{\mathbb{E}} \left[ L_{t}^{\tt charg}  \right] + V_{m}^{\tt charg} \mathop{\mathbb{E}} \left[k_{t} g(\mathbf{x}) \right]  \leq  D + c_{t}^{*} V_{m}^{\tt charg},\!
\end{equation}
where $D$ is the constant in (\ref{eq:drift_with_optimal}), and the optimal solution of $\mathbf{x}$ makes the term related to the virtual queues be $0$.
The above inequality holds for all $t>0$.
Therefore, we sum the above equation from $t=0$ to $t={\cal T}$, and then we can obtain
\begin{align} 
& \mathop{\mathbb{E}} \left[ L_{ {\cal T} +1 }^{\tt charg} \right]  - \mathop{\mathbb{E}} \left[ L_{0}^{\tt charg}  \right]  +  \sum_{t = 1}^{ \cal T } V_{m}^{\tt charg}  \mathop{\mathbb{E}} \left[k_{t} g(\mathbf{x}) \right]   \notag \\
& \leq  D {\cal T} + \sum_{t = 1}^{ \cal T } V_{m}^{\tt charg} c_{t}^{*}.
\end{align} 
As $\mathop{\mathbb{E}} [ L_{0}^{\tt charg} ]= 0$ and $\mathop{\mathbb{E}} [ L_{ {\cal T} +1 }^{\tt charg} ] > 0$ , we divide $V_{m}^{\tt charg}$ on both side and then we get the second inequality listed in Theorem \ref{theorem:optimal_gap}.

\end{document}